\documentclass[12pt,headings=normal]{scrartcl}

\usepackage{hyperref}
\usepackage[utf8]{inputenc}
\usepackage[ngerman]{babel}
\usepackage{amsmath,amssymb}
\usepackage[thmmarks]{ntheorem}
\usepackage{verbatim}
\usepackage{graphicx}
\usepackage[all]{xy}

\newcommand{\Ass}{\operatorname{Ass}}
\newcommand{\Ann}{\operatorname{Ann}}

\renewcommand{\P}{\mathcal{P}}
\renewcommand{\S}{\mathcal{S}}

\theorembodyfont{\itshape}
\theoremstyle{plain}
\newtheorem{Satz}{Satz}[section]
\newtheorem{Lemma}[Satz]{Lemma}
\newtheorem{Folgerung}[Satz]{Folgerung}

\theorembodyfont{\upshape}
\newtheorem{Bemerkung}[Satz]{Bemerkung}

\theoremheaderfont{\scshape}
\theorembodyfont{\upshape}
\theoremstyle{nonumberplain}
\theoremseparator{}
\theoremsymbol{$\Box$}
\newtheorem{Beweis}{Beweis}

\title{\Large \"Uber die von einem Ideal $I \subset R$ erzeugten $R$-Moduln II}
\author{\large Helmut Zöschinger\\
  \large Mathematisches Institut der Universität München\\
  \large Theresienstr. 39, D-80333 München\\
  \large E-mail: zoeschinger$@$mathematik.uni-muenchen.de
}
\date{}

\newcounter{abccount}
\newenvironment{abc}{%
\begin{list}{(\alph{abccount})}{%
\usecounter{abccount}%
\setlength{\partopsep}{0pt}%
\setlength{\topsep}{1ex}%
\setlength{\itemsep}{0pt}%
}%
}{\end{list}}

\newcounter{myenumcount}

\newcounter{iiicount}
\newenvironment{iii}{%
\begin{list}{(\roman{iiicount})}{%
\usecounter{iiicount}%
\setlength{\labelwidth}{3em}%
\setlength{\partopsep}{0pt}%
\setlength{\topsep}{1ex}%
\setlength{\itemsep}{0pt}%
}%
}{\end{list}}

\begin{document}
\maketitle

\centerline{\textbf{Abstract}}
\begin{abstract}
  \noindent
  Let $(R, \mathfrak m)$ be a commutative noetherian local ring and $I$ an
  ideal of $R$. Let $\P$ be the class of all $I$-generated $R$-modules $M$
  (i.e. there is an epimorphism $I^{(\Lambda)} \twoheadrightarrow M$) and let
  $\S$ be the class of all $I^{\circ}$-cogenerated $R$-modules $N$ (i.e.
  there is a monomorphism $N \hookrightarrow (I^{\circ})^{\Lambda}$ with
  $I^{\circ} = \operatorname{Hom}_R(I,E)$). We give a complete description
  of all injective and flat modules in $\P$ and $\S$. We show that
  $(\S,\P)$ forms a dual pair in the sense of Mehdi--Prest~(2015) and that
  $\P$ is always closed under pure submodules. We determine all ideals $I$
  for which $\P$ is closed under submodules, $\S$ is closed under factor
  modules and $\P$ (resp. $\S$) is closed under group extensions. In the last
  section, we examine the submodules $\gamma(M) = \sum\{U \subset M \,|\, U
  \in \P\}$ and $\kappa(M) = \bigcap \{V \subset M \,|\, M/V \in \S\}$ for all 
  $R$-modules $M$, and we specify their explicit structure in special cases. 
\end{abstract}

\bigskip

\noindent
\emph{Key words:} $I$-generated and $I^{\circ}$-cogenerated modules,
basically full ideals, dual pairs of modules, Matlis duality.
\bigskip

\noindent
\emph{Mathematics Subject Classification (2010):} 13C05, 13C11, 16D70, 16S90.

\section{Die injektiven $R$-Moduln in $\P$ und $\S$}
Stets sei $(R, \mathfrak{m})$ ein kommutativer, noetherscher, lokaler Ring,
$E$ die injektive Hülle des Restklassenkörpers $k =
R/\mathfrak{m}$ und $M^{\circ} = \operatorname{Hom}_R(M,E)$ das Matlis-Duale
eines $R$-Moduls $M$.

\begin{Lemma}\label{1.1}\mbox{}\\[-1em]
  \begin{abc}
  \item Für jeden $R$-Modul $M$ gelten die Implikationen
    \begin{equation*}
      IM=M \implies M \in \P \implies \Ann_R(I) \cdot M = 0.
    \end{equation*}
  \item Ist $M$ \emph{injektiv} und $\overline{I} = \Ann_R\Ann_R(I)$, so ist
    jede der drei Bedingungen äquivalent mit $M[\overline{I}] = 0$.
  \item Für jedes Primideal $\mathfrak{p}$ von $R$ gilt
    \begin{equation*}
      \overline{I} \not\subset \mathfrak{p} \iff \Ann_R(I) \cdot
      R_{\mathfrak{p}} = 0 \iff I_{\mathfrak{p}} \ \text{ist regulär
        in}\ R_{\mathfrak{p}}.
    \end{equation*}
  \end{abc}
\end{Lemma}

\begin{Beweis}
  (a) ist klar, weil nach (\cite{007} Proposition~1.1) $M$ genau dann zu
  $\P$ gehört, wenn es eine Erweiterung $M \subset X$ gibt mit $M =IX$.\\
  (b) Für jeden $R$-Modul $M$ und jedes Ideal $\mathfrak{a}$ von $R$ sei
  $M[\mathfrak{a}] := \Ann_M(\mathfrak{a}) = \{x \in M \,|\, ax = 0 \
  \text{für alle}\ a \in \mathfrak{a}\}$. Ist nun $M$ injektiv, gilt bekanntlich
  \begin{equation*}
    \mathfrak{a} M = M[\Ann_R(\mathfrak{a})],
  \end{equation*}
  so dass aus $\Ann_R(I) \cdot M = 0$, d.\,h. $M = M[\Ann_R(I)]$ sofort
  folgt $IM=M$, wegen $\Ann_R(I) \cdot M = M[\overline{I}]$ aber auch der
  zweite Teil.\\
  (c) Beide Äquivalenzen sind wohlbekannt: Genau dann ist
  $\mathfrak{a}R_{\mathfrak{p}} = 0$, wenn es ein $s \in R\setminus
  \mathfrak{p}$ gibt mit $s \mathfrak{a} = 0$, und $\Ann_R(I) \cdot
  R_{\mathfrak{p}}$ ist im Ring $R_{\mathfrak{p}}$ der Annullator des Ideals
  $I R_{\mathfrak{p}}$.
\end{Beweis}

Für jeden $R$-Modul $X$ und jedes Ideal $\mathfrak{a}$ von $R$ ist genau
dann $X[\mathfrak{a}] = 0$, wenn $\mathfrak{a} \not\subset \mathfrak{p}$ ist
für alle $\mathfrak{p} \in \Ass(X)$, so dass mit (b) und (c) folgt:

\begin{Satz}\label{1.2}
  Ist $M$ injektiv, so gilt: $M \in \P \iff I_{\mathfrak{p}}$ ist regulär in
  $R_{\mathfrak{p}}$ für alle $\mathfrak{p} \in \Ass(M)$.
\end{Satz}

\begin{Bemerkung}\label{1.3}
  In (1.1, b) kann man nicht erwarten, dass $M[I] = 0$ ist. Ist z.\,B. $V$
  der kleinste Untermodul von $E$ mit $E/V \in \S$, so gilt $E[\overline{I}]
  \subset V \subset E[I]$, und der erste Untermodul ist genau dann Null,
  wenn $I$ regulär ist, der zweite, wenn $I \cong R$ ist (\cite{007}~p.~6,
  Beispiel~4), der dritte, wenn $I=R$ ist.\hfill$\Box$
\end{Bemerkung}

\begin{Bemerkung}\label{1.4}
  Für jeden \emph{flachen} $R$-Modul $N$ ist $N^{\circ} =
  \operatorname{Hom}_R(N,E)$ injektiv, so dass man alle bisherigen
  Ergebnisse via Matlis-Dualität von $N^{\circ}$ auf $N$ übertragen kann.
  Weil nach (\cite{007} Proposition~3.1) $N^{\circ} \in \P$ äquivalent ist
  mit $N \in \S$, erhält man sofort:\\
  Für jeden $R$-Modul $N$ gilt: $N[I] = 0 \implies N \in \S \implies
  \Ann_R(I) \cdot N = 0$. Falls $N$ \emph{flach} ist, ist jede der drei
  Bedingungen äquivalent mit $\overline{I} \cdot N = N$, d.\,h. mit
  $I_{\mathfrak{p}}$ regulär in $R_{\mathfrak{p}}$ für alle $\mathfrak{p}
  \in \operatorname{Koass}(N)$.\hfill$\Box$
\end{Bemerkung}

Ist $N$ injektiv und $N \in \S$, folgt aus $\Ann_R(I) \cdot N = 0$ nach
(1.1, b) $N \in \P$. Um wie viel stärker die Bedingung $N \in \S$ ist,
wollen wir jetzt präzisieren:

\begin{Lemma}\label{1.5}
  Für ein Primideal $\mathfrak{p}$ von $R$ sind äquivalent:
  \begin{iii}
    \item $E(R/\mathfrak{p}) \in \S$
    \item $(R_{\mathfrak{p}})^{\circ\circ} \in \P$
    \item $R_{\mathfrak{p}} \in \P$
    \item $I_{\mathfrak{p}} \cong R_{\mathfrak{p}}$
  \end{iii}
\end{Lemma}

\begin{Beweis}
  (i $\to$ iv) Nach Schenzel (\cite{005} Lemma~2.3) ist
  $E(R/\mathfrak{p})^{\circ} \cong \widehat{R_{\mathfrak{p}}^{(J)}}$ für eine
  Indexmenge $J \neq \varnothing$, also
  \begin{equation*}
    \widehat{R_{\mathfrak{p}}}\, \subset^{\oplus}\, E(R/\mathfrak{p})^{\circ},
  \end{equation*}
  wobei $\widehat{\phantom{aa}}$ die Vervollständigung über
  $R_{\mathfrak{p}}$ sei. Mit $E(R/\mathfrak{p})^{\circ}$ ist also auch $X =
  \widehat{R_{\mathfrak{p}}}$ aus $\P$, der $R$-Epimorphismus $I^{(\Lambda)}
  \twoheadrightarrow X$ induziert einen $\widehat{R}$-Epimorphismus
  $(I_{\mathfrak{p}})^{(\Lambda)} \twoheadrightarrow X$, und weil $X$ über
  $R_{\mathfrak{p}}$ das Biduale von $R_{\mathfrak{p}}$ ist, folgt nach (\cite{007}
  p.~6, Beispiel~4) $I_{\mathfrak{p}} \cong R_{\mathfrak{p}}$.\\
  (iv $\to$ iii) Es genügt zu zeigen, dass $I_{\mathfrak{p}}$ als $R$-Modul
  $I$-generiert ist. Allgemeiner gilt aber für jeden $R$-Modul $A \in \P$,
  d.\,h. $A \cong IB$, dass $A_S \cong I \cdot B_S$, also auch
  $A_S \in \P$ ist.\\
  (iii $\to$ ii) Für jeden $R$-Modul $A \in \P$, d.\,h. $A \cong IB$, ist nach
  (\cite{007} Proposition~3.1) $A^{\circ\circ} \cong I\cdot B^{\circ\circ}$,
  also auch $A^{\circ\circ} \in \P$.\\
  (ii $\to$ i) Es genügt zu zeigen, dass
  \begin{equation*}
    E(R/\mathfrak{p}) \,\subset^{\oplus}\, (R_{\mathfrak{p}})^{\circ}
  \end{equation*}
  ist. Bei $\mathfrak{p} = \mathfrak{m}$ ist sogar $E \cong R^{\circ}$, bei
  $\mathfrak{p} \neq \mathfrak{m}$ gilt nach (\cite{006} p.~7, l.~2)\\
  $(R_{\mathfrak{p}})^{\circ} \cong \coprod\limits_{\mathfrak{q} \subset
    \mathfrak{p}} E(R/\mathfrak{q})^{(J_{\mathfrak{q}})}$ mit
  $J_{\mathfrak{q}} \neq \varnothing$ für alle $\mathfrak{q} \subset
  \mathfrak{p}$, also wieder die Behauptung.
\end{Beweis}

Weil jeder injektive $R$-Modul direkte Summe von unzerlegbaren der Form
$E(R/\mathfrak{p})$ ist, erhält man:

\begin{Satz}\label{1.6}
  Ist $N$ injektiv, so gilt: $N \in \S \iff I_{\mathfrak{p}} \cong
  R_{\mathfrak{p}}$ für alle $\mathfrak{p} \in \Ass(N)$.
\end{Satz}

\begin{Bemerkung}
  Die Implikation (ii $\to$ iii) in (\ref{1.5}) gilt sogar für jeden
  $R$-Modul $M$, denn nach (\ref{2.1}, a) ist $M \in \P$ äquivalent mit
  $M^{\circ} \in \S$. Falls also $M$ \emph{flach} war, ist $M^{\circ}$
  injektiv, und (\ref{1.6}) liefert sofort: $M \in \P \iff I_{\mathfrak{p}}
  \cong R_{\mathfrak{p}}$ für alle $\mathfrak{p} \in
  \operatorname{Koass}(M)$.\hfill$\Box$
\end{Bemerkung}

\section{Das duale Paar $(\S, \P)$}

Die Klasse $\P$ aller $I$-generierten $R$-Moduln ist natürlich gegenüber
direkten Produkten, direkten Summen und Faktormoduln abgeschlossen. Weiter
gilt:

\begin{Lemma}\label{2.1}
  Die Klasse $\P$ ist gegenüber
  \begin{abc}
  \item reinen Untermoduln,
  \item reinen Gruppenerweiterungen und
  \item rein-injektiven Hüllen
  \end{abc}
  abgeschlossen.
\end{Lemma}

\begin{Beweis}
  (a) \emph{1.~Schritt} Eine beliebige Modulerweiterung $M \subset X$ heißt
  nach (\cite{007} p.~7) \emph{$I$-klein}, wenn $M_{*} = I(M :_X I)$ mit $M$
  übereinstimmt. War $X$ \emph{injektiv}, ist das nach (\cite{007}
  Proposition~1.1, iii) äquivalent mit $M \in \P$.\\
  Dual heißt $B \subset A$ nach (\cite{007} p.~6) \emph{$I$-groß}, wenn $B$
  mit $B^{*} = (IB) :_A I$ übereinstimmt (siehe auch \cite{002} Theorem
  2.12). Allgemeiner als in \cite{007} gilt jetzt: War $A$ \emph{flach}, ist
  $B = B^{*}$ äquivalent mit $A/B \in \S$. Nur "`$\Leftarrow$"' ist zu
  zeigen, und da ist in der exakten Folge
  \begin{equation*}
    0 \longrightarrow (A/B)^{\circ} \longrightarrow A^{\circ}
    \longrightarrow B^{\circ} \longrightarrow 0
  \end{equation*}
  $(A/B)^{\circ} \in \P$ und $A^{\circ}$ injektiv, also nach eben
  $\Ann_{A^{\circ}}(B)$ $I$-klein in $A^{\circ}$ und dann nach (\cite{007}
  Bemerkung~4.4) $B$
  $I$-groß in $A$.\\
  \emph{2.~Schritt} Zeigen wir zuerst den Spezialfall $M \subset
  M^{\circ\circ}$ mit $M^{\circ\circ} \in \P$, d.\,h. $M^{\circ} \in \S$.
  Mit einer injektiven Erweiterung $M \subset X$ ist dann
  \begin{equation*}
    0 \longrightarrow (X/M)^{\circ} \longrightarrow X^{\circ}
    \longrightarrow M^{\circ} \longrightarrow 0
  \end{equation*}
  exakt und $X^{\circ}$ flach, also nach eben $\Ann_{X^{\circ}}(M)$ $I$-groß
  in $X^{\circ}$, und mit (\cite{007} Lemma~4.3) folgt $M$ $I$-klein in $X$,
  also $M \in \P$.\\
  Ist jetzt $A \hookrightarrow B$ ein reiner Monomorphismus und $B \in \P$,
  ist mit $B^{\circ\circ} \in \P$ auch der direkte Summand $A^{\circ\circ}
  \in \P$, also $A \in \P$.\\
  (b) Ist $0 \to A \to B \to C \to 0$ rein-exakt und $A,C \in \P$, folgt aus
  $C^{\circ}, A^{\circ} \in \S$ sogar $B^{\circ} \in \S$ (weil die
  entsprechende Folge zerfällt), also nach (a) $B \in \P$.\\
  (c) Sei $M \in \P$ und $M \subset N$ eine rein-injektive Hülle. Dann ist
  $N$ bis auf Isomorphie direkter Summand in $M^{\circ\circ}$, also auch $N
  \in \P$.
\end{Beweis}

Aus (\ref{2.1}) und den Regeln
\begin{equation*}
  N \in \S \iff N^{\circ} \in \P, \quad M \in \P \iff M^{\circ} \in S
\end{equation*}
folgt, dass $(\S,\P)$ ein im Sinne von Mehdi--Prest (\cite{003} p.~1389)
\emph{duales Paar} bildet. Wir wollen untersuchen, wann seine Komponenten
gegenüber Untermoduln, Faktormoduln oder Gruppenerweiterungen abgeschlossen
sind.

\begin{Satz}\label{2.2}
  Äquivalent sind:
  \begin{iii}
    \item $\P$ ist gegenüber Untermoduln abgeschlossen.
    \item $\S$ ist gegenüber Faktormoduln abgeschlossen.
    \item $\P = \S = \{M \in R\text{-}\mathrm{Mod}\ |\, \Ann_R(I) \cdot M =
      0\}$.
    \item Es gibt einen Epimorphismus $I \twoheadrightarrow R/\Ann_R(I)$.
  \end{iii}
\end{Satz}

\begin{Beweis}
  (i $\to$ iv) Mit $I = R r_1 + \dots + R r_n$ ($n \geq 1$) und $x =
  (r_1,\dotsc, r_n) \in I \times \dots \times I$ ist nach Voraussetzung auch
  $Rx \in \P$, d.\,h. $I \twoheadrightarrow Rx \cong R/\Ann_R(x) =
  R/\Ann_R(I)$.\\
  (iv $\to$ iii) In $\P \subset \{M \in R\text{-}\mathrm{Mod}\ |\, \Ann_R(I)
  \cdot M = 0\}$ gilt Gleichheit, denn aus $R^{(\Lambda)} \twoheadrightarrow
  M$ und $\Ann_R(I) \cdot M = 0$ folgt $(R/\Ann_R(I))^{(\Lambda)}
  \twoheadrightarrow M$, also nach Voraussetzung $I^{(\Lambda)}
  \twoheadrightarrow M$. Damit gilt auch $N \in \S \iff N^{\circ} \in \P
  \iff \Ann_R(I) \cdot N^{\circ} = 0 \iff \Ann_R(I) \cdot N = 0$.\\
  (iii $\to$ i) klar.\\
  (i $\to$ ii) Aus $N \in \S$ und $V \subset N$ folgt $N^{\circ} \in \P$,
  also nach Voraussetzung auch $(N/V)^{\circ} \in \P$, d.\,h. $N/V \in
  \S$.\\
  (ii $\to$ iv) Mit $I$ und $x$ wie am Anfang folgt $(I^{\circ})^n
  \twoheadrightarrow (Rx)^{\circ}$, so dass nach Voraussetzung $(Rx)^{\circ}
  \in \S$, d.\,h. $Rx \in \P$ ist, also $I \twoheadrightarrow R/\Ann_R(I)$.
\end{Beweis}

\begin{Bemerkung}
  Die Bedingung (iv) ist natürlich erfüllt, wenn $I$ zyklisch oder
  halbeinfach $\neq 0$ ist. Das ist z.\,B. der Fall, wenn $\mathfrak{m}^2 =
  0$ ist.\hfill$\Box$
\end{Bemerkung}

\begin{Bemerkung}
  Die Bedingung (iv) ist genau dann erfüllt, wenn es einen \emph{zyklischen}
  direkten Summanden $I_1$ von $I$ gibt mit $\Ann_R(I_1) = \Ann_R(I)$. War
  also $I$ direkt unzerlegbar (z.\,B. wenn $R$ uniform ist), muss $I$
  bereits zyklisch sein.\hfill$\Box$
\end{Bemerkung}

\begin{Satz}\label{2.5}
  Äquivalent sind:
  \begin{iii}
    \item $\P$ ist gegenüber Gruppenerweiterungen abgeschlossen.
    \item $\S$ ist gegenüber Gruppenerweiterungen abgeschlossen.
    \item Es ist $I = 0$ oder $I \cong R$.
  \end{iii}
\end{Satz}

\begin{Beweis}
  (iii $\to$ i) ist klar, ebenso (i $\to$ ii), denn aus $0 \to A \to B \to C
  \to 0$ exakt und $A,C \in \S$ folgt $C^{\circ}, A^{\circ} \in \P$, also
  nach Voraussetzung $B^{\circ} \in \P$, d.\,h. $B \in \S$.\\
  Für (ii $\to$ iii) sei gleich $I \neq 0$. Zeigen wir im \emph{1.~Schritt},
  dass $I$ regulär ist: Mit $R/\mathfrak{m}$ gehört auch jeder reduzierte
  $R$-Modul $N$ zu $\S$, denn der kleinste Untermodul $V$ von $N$, mit $N/V
  \in \S$, hat nach Voraussetzung keinen maximalen Untermodul, ist also
  radikalvoll, d.\,h. es ist $V = 0$. Insbesondere ist $R \in \S$, also
  $\Ann_R(I) = 0$.\\
  Sei im \emph{2.~Schritt} $V$ der kleinste Untermodul von $E$ mit $E/V \in
  \S$. Dann ist $\operatorname{Hom}_R(V,I^{\circ}) = 0$, denn für jedes $f
  \colon V \to I^{\circ}$ ist $V/\operatorname{Ke} f \in \S$, also nach
  Voraussetzung $E/\operatorname{Ke} f \in \S \implies \operatorname{Ke} f =
  V$, d.\,h. $f=0$. Aus $\operatorname{Hom}_R(I,V^{\circ}) = 0$ folgt
  $\Ass(V^{\circ}) \cap \operatorname{Supp}(I) = \varnothing$ (\cite{001}
  chap.~IV, \verb+§+~1, Proposition~10), wegen $\operatorname{Supp}(I) =
  V(\Ann_R(I)) = \operatorname{Spec}(R)$ nach dem 1.~Schritt also
  $\Ass(V^{\circ}) = \varnothing$, $V=0$. Aber $E \in \S$ bedeutet nach
  (\cite{007} p.~6, Beispiel~4) schon $I \cong R$.
\end{Beweis}

\section{Die Funktoren $\gamma$ und $\kappa$}

Für jeden $R$-Modul $M$ hat die Menge $\{U \subset M \,|\, U \in \P\}$ ein
größtes Element, das wir mit $\gamma(M)$ bezeichnen, und $\{V \subset M
\,|\, M/V \in \S\}$ ein kleinstes Element, das wir mit $\kappa(M)$
bezeichnen (siehe die beiden Schritte im letzten Beweis).

\medskip

\noindent
\textbf{Beispiel 1} Gibt es einen Epimorphismus $I \twoheadrightarrow
R/\Ann_R(I)$ wie in (\ref{2.2}), so gilt für jeden $R$-Modul $M$
\begin{equation*}
  \gamma(M) = M[\Ann_R(I)], \quad \kappa(M) = \Ann_R(I) \cdot M.
\end{equation*}

\begin{Beweis}
  Nach (\ref{2.2}, iii) gilt für $U \subset M$, dass $U \in \P$ äquivalent
  ist mit $U \subset M[\Ann_R(I)]$, und für $V \subset M$, dass $M/V \in \S$
  äquivalent ist mit $\Ann_R(I) \cdot M \subset V$.
\end{Beweis}

\noindent
\textbf{Beispiel 2}
\begin{abc}
\item Ist $M \subset X$ und $X$ injektiv, so gilt $\gamma(M) = I(M :_X I)$.
\item Ist $B \subset A$ und $A$ flach, so gilt $\kappa(A/B) = (IB :_A I)/B$.
\end{abc}

\begin{Beweis}
  (a) "`$\subset$"' Für $U = \gamma(M)$ gilt nach (\cite{007}
  Proposition~1.1, iii) $U = I(U :_X I) \subset I(M :_X I) = M_{*}$.
  "`$\supset$"' $M_{*}$ ist ein Untermodul von $M$ und natürlich
  $I$-generiert, also enthalten in $\gamma(M)$.\\
  (b) "`$\supset$"' Mit $V/B  = \kappa(A/B)$ ist $A/V \in \S$, also nach dem
  Beweis von (\ref{2.1}, a) $V = (IV) :_A I \supset (IB) :_A I = B^{*}$.
  "`$\subset$"' $\frac{A}{IB} / \frac{A}{IB}[I]$ ist nach (\cite{007}
  Proposition~3.1) $I^{\circ}$-kogeneriert und isomorph zu $A/B^{*} \cong
  \frac{A}{B} / \frac{B^{*}}{B^{\phantom{*}}}$, also $\kappa(\frac{A}{B}) \subset
  \frac{B^{*}}{B^{\phantom{*}}}$.
\end{Beweis}

\noindent
\textbf{Beispiel 3} Ist $R$ ein Integritätsring und
$\operatorname{Ext}_R^1(R/I, R) = 0$, so gilt $\gamma(R) = I$.

\begin{Beweis}
  Bei $I = 0$ ist $\gamma(M) = 0$ für jeden $R$-Modul $M$. Bei $I \neq 0$
  folgt mit dem Quotientenkörper $K$ von $R$, dass in der exakten Folge
  \begin{equation*}
    \operatorname{Hom}_R(R/I, K) \longrightarrow \operatorname{Hom}_R(R/I,
    K/R) \longrightarrow \operatorname{Ext}_R^1(R/I, R)
  \end{equation*}
  das erste und dritte Glied Null ist, also auch $(K/R)[I] = (R :_K I)/R$.
  Aus $R :_K I = R$ folgt dann mit Beispiel~2~(a) die Behauptung.
\end{Beweis}

\begin{Satz}\label{3.1}
  Für jeden $R$-Modul $M$ gilt:
  \begin{abc}
  \item $IM \subset \gamma(M) \subset M[\Ann_R(I)]$ und $\Ann_R(I) \cdot M
    \subset \kappa(M) \subset M[I]$.
  \item $\gamma(M)$ ist groß in $M[\Ann_R(I)]$ und $\kappa(M)/\Ann_R(I)
    \cdot M$ ist klein in $M/\Ann_R(I) \cdot M$.
  \item $\Ass(\gamma(M)) = \Ass(M) \cap \operatorname{Supp}(I)$ und
    $\operatorname{Koass}(M/\kappa(M)) = \operatorname{Koass}(M) \cap
    \operatorname{Supp}(I)$.
  \end{abc}
\end{Satz}

\begin{Beweis}
  (a) Aus $IM \in \P$ und $\gamma(M) \in \P$ folgen die beiden ersten
  Inklusionen, aus $M/\kappa(M) \in \S$ und $M/M[I] \in \S$ (\cite{007}
  Proposition~3.1) die beiden anderen.\\
  (b) Sei $0 \neq U \subset M[\Ann_R(I)]$. Mit irgendeinem $\mathfrak{p}_0 \in
  \Ass(U)$ ist dann $\Ann_R(I) \cdot R/\mathfrak{p}_0 = 0$, also
  $\mathfrak{p}_0 \in \Ass(U) \cap \operatorname{Supp}(I) =
  \Ass(\operatorname{Hom}_R(I,U))$ nach (\cite{001} chap.~IV, \verb+§+~1,
  Proposition~10) $\implies \operatorname{Hom}_R(I,U) \neq 0$, $\gamma(U)
  \neq 0$.\\ Entsprechend müssen wir für jeden Zwischenmodul $\Ann_R(I)\cdot M
  \subset V \subsetneqq M$ zeigen, dass $V + \kappa(M) \subsetneqq M$ ist.
  Mit irgendeinem $\mathfrak{p}_0 \in \operatorname{Koass}(M/V)$ folgt, weil
  $M/V$ durch $\Ann_R(I)$ annulliert wird, $\mathfrak{p}_0 \in
  \Ass((M/V)^{\circ}) \cap \operatorname{Supp}(I)$, also
  $\operatorname{Hom}_R(I, (M/V)^{\circ}) \cong \operatorname{Hom}_R(M/V,
  I^{\circ}) \neq 0$, $V \subset V_1 \subsetneqq M$ mit $M/V_1 \in \S$, also
  $\kappa(M) \subset V_1$.\\
  (c) Für jeden wesentlichen Monomorphismus $A \hookrightarrow B$ gilt
  bekanntlich $\Ass(A) = \Ass(B)$, für jeden wesentlichen Epimorphismus $B
  \twoheadrightarrow C$ dual $\operatorname{Koass}(B) =
  \operatorname{Koass}(C)$.
\end{Beweis}

\begin{Folgerung}\label{3.2}
  Mit $\overline{I} = \Ann_R \Ann_R(I)$ gilt:
  \begin{abc}
  \item $M$ injektiv $\implies IM = \gamma(M) = M[\Ann_R(I)]$ und
    $M[\overline{I}] \subset \kappa(M) \subset M[I]$.
  \item $M$ flach $\implies \Ann_R(I) \cdot M = \kappa(M) = M[I]$ und $IM
    \subset \gamma(M) \subset \overline{I}M$.
  \end{abc}
\end{Folgerung}

\begin{Beweis}
  (a) Weil $M$ injektiv ist, gilt $\mathfrak{a} M = M[\Ann_R(\mathfrak{a})]$
  für alle Ideale $\mathfrak{a}$ von $R$, so dass $\mathfrak{a} = I$ in
  (\ref{3.1}, a) die Gleichung und $\mathfrak{a} = \Ann_{R}(I)$ die
  Ungleichung liefert.\\
  (b) Entsprechend, weil für jeden flachen $R$-Modul $M$ gilt
  $\Ann_R(\mathfrak{a}) \cdot M = M[\mathfrak{a}]$.
\end{Beweis}

\begin{Folgerung}\label{3.3}
  Falls $I$ ein Annullatorideal, d.\,h. $I = \overline{I}$ ist, folgt
  \begin{abc}
  \item für jeden injektiven $R$-Modul $M \in \P$ sogar $M \in \S$,
  \item für jeden flachen $R$-Modul $M \in \S$ sogar $M \in \P$.
  \end{abc}
\end{Folgerung}

\begin{Beweis}
  Bei (a) ist $M[\overline{I}] =0$ nach (\ref{1.1}, b), also $\kappa(M)=0$,
  bei (b) ist $\overline{I} \cdot M = M$ nach (\ref{1.4}), also $\gamma(M) =
  M$.
\end{Beweis}

\begin{Bemerkung}
  (1) Auch wenn $I$ zyklisch ist, kann $IM \subsetneqq \gamma(M)$ und
  $\kappa(M) \subsetneqq M[I]$ sein, z.\,B. wenn $I \subset  \Ann_R(M)
  \subsetneqq R$ und $I \cong R$ ist (denn dann ist $IM=0$, $\gamma(M)=M$,
  $\kappa(M) = 0$ und $M[I] = M$).\\
  (2) Ist $I$ nicht zyklisch, kann auch $\gamma(M) \subsetneqq M[\Ann_R(I)]$
  und $\Ann_{R}(I) \cdot M \subsetneqq \kappa(M)$ sein, z.\,B. wenn $I$
  regulär, aber nicht isomorph zu $R$ ist (denn dann ist $\gamma(R)
  \subsetneqq R$ und $0 \subsetneqq \kappa(E)$).\hfill$\Box$
\end{Bemerkung}

Ist ein $R$-Modul $M$ einreihig (d.\,h. $\mathcal{L}(M)$ totalgeordnet und
endlich), so zeigten wir in (\cite{007} Satz~1.4 und Folgerung~3.3), dass
$M\in\P$ äquivalent ist mit $M \in \S$ und das (mit Länge$(M)
= n \geq 1$) weiter äquivalent ist mit $\mathfrak{m}^{n-1} \cdot I
\not\subset \Ann_R(M) \cdot I$. Wir wollen zum Schluss genauer $\gamma(M)$
und $\kappa(M)$ berechnen, ebenso (weil auch $M^{\circ}$ einreihig ist)
$\gamma(M^{\circ})$ und $\kappa(M^{\circ})$.

\begin{Satz}\label{3.5}
  Sei $M$ einreihig von der Länge $n \geq 1$, $M \supsetneqq M_1 \supsetneqq
  M_2 \supsetneqq \dots \supsetneqq M_{n-1} \supsetneqq 0$ sein
  Untermodulverband und $s$ die kleinste natürliche Zahl mit $\mathfrak{m}^s
  \cdot I \subset \Ann_R(M) \cdot I$. Dann ist $0 \leq s \leq n$ und
  \begin{equation*}
    \gamma(M) = M_{n-s},\quad \kappa(M) = M_s.
  \end{equation*}
\end{Satz}

\begin{Beweis}
  Aus $\mathfrak{m}^n \subset \Ann_R(M)$ folgt $s \leq n$. Nach (\cite{007}
  Folgerung~1.5) gilt für alle $0\leq i \leq n$, dass $M/M_i \notin \P$
  äquivalent ist mit $i > 0$ und $\mathfrak{m}^{i-1} \cdot I \subset
  \Ann_R(M) \cdot I$, d.\,h. mit $i-1 \geq s$, $i > s \implies \kappa(M) =
  M_s$.\\
  Weil auch $M^{\circ}$ einreihig von der Länge $n \geq 1$ ist und
  $(M_i)^{\circ}$ ein Faktormodul der Länge $n-i$, ist nach demselben Zitat
  $M_i \notin \P$ äquivalent mit $n-i > 0$ und $\mathfrak{m}^{(n-i)-1} \cdot I
  \subset \Ann_R(M) \cdot I$, d.\,h. mit $(n-i)-1 \geq s$, $n-s > i \implies
  \gamma(M) = M_{n-s}$.
\end{Beweis}

\begin{Folgerung}
  Sei $M$ einreihig von der Länge $n \geq 1$. Dann ist
  \begin{equation*}
    \gamma(M^{\circ}) = \Ann_{M^{\circ}}(\kappa(M)), \quad \kappa(M^{\circ})
    = \Ann_{M^{\circ}}(\gamma(M)).
  \end{equation*}
\end{Folgerung}

\begin{Beweis}
  Auch $D = M^{\circ}$ ist einreihig von der Länge $n \geq 1$, genauer $D
  \supsetneqq D_1 \supsetneqq D_2 \supsetneqq \dots \supsetneqq D_{n-1}
  \supsetneqq 0$ sein Untermodulverband mit $D_i =
  \Ann_{M^{\circ}}(M_{n-i})$. Wegen $\Ann_R(D) = \Ann_R(M)$ folgt mit
  (\ref{3.5}) sofort $\gamma(M^{\circ}) = D_{n-s} = \Ann_{M^{\circ}}(M_s) =
  \Ann_{M^{\circ}}(\kappa(M))$, entsprechend
  $\kappa(M^{\circ})$.
\end{Beweis}


\begin{thebibliography}{99}
\bibitem{001}
  N. Bourbaki: \emph{Algèbre commutative}: Hermann. Paris (1967)
\bibitem{002}
  W. J. Heinzer -- L. J. Ratliff Jr. -- D. E. Rush: \emph{Basically full
    ideals in local rings}: J.~Algebra 250 (2002) 371--396
\bibitem{003}
  A. R. Mehdi -- M. Prest: \emph{Almost dual pairs and definable classes of
    modules}: Commun. Algebra 43 (2015) 1387--1397
\bibitem{004}
  J. J. Rotman: \emph{An introduction to homological algebra}: Academic
  Press. New York (1997)
\bibitem{005} P. Schenzel: \emph{A note on the Matlis dual of a certain
    injective hull}: J. pure appl. Algebra 219 (2015) 666--671  
\bibitem{006}
   H. Zöschinger: \emph{Über rein-wesentliche Erweiterungen}: arXiv
   1403.5957 (2014) 1--11
\bibitem{007}  H. Z"oschinger: \emph{Über die von einem Ideal $I \subset R$
    erzeugten $R$-Moduln}: arXiv 1604.02349 (2016) 1--9
\end{thebibliography}
\end{document}